\documentclass{amsart}
\usepackage{amsmath}
\usepackage{amsfonts}
\usepackage{amssymb}
\usepackage{amsbsy,amsthm,latexsym,
            amsopn,amstext,amsxtra,euscript,amscd,amsthm}

\newtheorem{proposition}{\textbf{Proposition}}
\newtheorem{theorem}{\textbf{Theorem}}
\newtheorem{corollary}{\textbf{Corollary}}
\newtheorem{remark}{\textbf{Remark}}

\newtheorem{lemma}{\textbf{Lemma}}
\newtheorem{question}{\textbf{Question}}
\newtheorem{exe}{\textbf{Example}}

\def\N {\mathbb{N}}
\def\Z {\mathbb{Z}}
\def\Q {\mathbb{Q}}

\def\R {\mathbb{R}}

\def\cR {\mathcal{R}}

\def\QQ {\overline{\Q}}
\def\C {\mathbb{C}}
\def\cF {\mathcal{F}}

\theoremstyle{remark}

\numberwithin{equation}{section}

\begin{document}

\title[FIELD THEORETIC PROPERTIES AND TRANSCENDENTAL NUMBERS]{SOME FIELD THEORETIC PROPERTIES AND AN APPLICATION CONCERNING TRANSCENDENTAL NUMBERS}

\author{CHRISTIAN U. JENSEN}
\address{DEPARTMENT OF MATHEMATICAL SCIENCES, UNIVERSITY OF COPENHAGEN,
COPENHAGEN, DENMARK}
\email{cujensen@math.ku.dk}
\author{DIEGO MARQUES}\thanks{The second author thanks to CAPES for the financial support.}
\address{DEPARTAMENTO DE MATEM\'{A}TICA, UNIVERSIDADE DE BRAS\' ILIA,
BRAS\' ILIA, DF, BRAZIL}
\email{diego@mat.unb.br}

\subjclass[2000]{12F10; 11J81}

\keywords{Artin-Schreier, Galois theory, Gelfond-Schneider, Transcendental numbers}

\begin{abstract}
For a proper subfield $K$ of $\QQ$ we show the
existence of an algebraic number $\alpha$ such that no power
$\alpha^n$, $n\geq 1$, lies in $K$. As an application it is shown that
these numbers, multiplied by convenient Gaussian numbers, can be written in the form $P(T)^{Q(T)}$ for some transcendental numbers $T$
where $P$ and $Q$ are arbitrarily prescribed non-constant rational functions over
$\QQ$.
\end{abstract}

\maketitle

\section{Introduction}
The origin of this paper was a result concerning the transcendency of
some numbers appearing  in certain exponential equations. This lead to
some purely field theoretic questions, which may be of some
independent interest. The first of the results are easy
consequences of the famous
\begin{theorem}\label{AS}[Artin-Schreier]
(Cf. \cite{Jac}) Let $M$ be an algebraically closed
field. If the characteristic of $M$ is zero, any proper subfield $F$
of $M$ such that $[M:F]$ is finite, is real closed (i.e., $F$  is not
algebraically closed but $F(\sqrt{-1})$ is algebraically closed) and thus $[M:F] =
2$. If the characteristic of $M$ is a prime number there is no proper
subfield $F$ of $M$ for which $[M:F]$ is finite. 
\end{theorem}

\begin{remark}
Artin-Schreier's theorem implies that the absolute Galois
group of any non real closed field contains no automorphisms of finite
order ($\neq 1$).
\end{remark}

In Section 2 we give some field theoretic results to be used in the final
Section 3. This section first brings a brief historical survey of
classical
famous theorems concerning transcendental numbers and then - as an
application of the results in Section 2 - we show the existence of
algebraic numbers that can be written as powers of two transcendental
numbers of a very special form.

\section{Some field theoretic properties}

We start with some consequences of the Theorem \ref{AS}. Although we shall only need fields of characteristic zero
in our later applications, for the sake of completeness, we consider also the finite characteristic case in the first three propositions.
\begin{proposition}\label{p1}
Let $M$ be an algebraically closed  field and $F$ a proper subfield which is not real closed. Then there exist elements in $M$ of arbitrarily large degree with respect to $F$.
\end{proposition}
\begin{proof}
In the case where the characteristic of $M$ (= the
characteristic of $F$) is a prime number $p$ we may assume that $F$
 is a perfect field. Indeed,
if there were an element $a\in F$ such that $a$ were not the $p$-th
power of an element in $F$ then the polynomials $x^{p^e} -a$ were irreducible in $F[x]$
for {\it every} natural number $e$. This would give rise to elements
in $M$ of arbitrary high degrees with respect to $F$.

By Theorem 1 $M$ is an infinite algebraic extension of $F$. Hence there exists
an infinite strictly increasing sequence of finite extensions
contained in $M$
 $$
 F \subsetneq M_1 \subsetneq M_2 \subsetneq \dots \subsetneq M_i\
\subsetneq M_{i +1} \subsetneq \dots $$

Since $F$ is perfect all these extensions are separable, so if $\alpha_i$ is a primitive
element for $M_i/F$ the degree of $\alpha_i$ with respect to $F$ tends
to infinity as $i$ tends to infinity.
\end{proof}
Since real closed fields necessarily have characteristic $0$,
Proposition 1 implies
\begin{proposition}
Let $F$ be an arbitrary field. If the characteristic
of $F$ is zero and there exists an irreducible
polynomial in $F[x]$ of degree $>2$, then there exist irreducible
polynomials in $F[x]$ of arbitrarily large degrees.  

If the characteristic of $F$ is a prime number and there exists an
irreducible polynomial in $F[x]$ of degree $> 1$, then there exist
irreducible polynomials in $F[x]$ of arbitrarily large degrees. 
\end{proposition}
\begin{proposition}\label{p2}
Let $M$ be an algebraically closed field and $F$ a proper subfield which is not real closed. Then $M$ cannot be generated by adjoining to $F$ all elements having degree (with respect to $F$) at most some fixed number $n$.
\end{proposition}
\begin{proof}
We may obviously assume that $M$ is an algebraic
extension of $F$.
The  Galois closure of the field obtained by
adjoining all elements of degree (w.r.t. $F$) at most $n$ has a
Galois group which is a subgroup of a direct product of symmetric
groups of degrees $\leq n$. This Galois closure cannot be $M$, since
otherwise the absolute Galois group of $F$ would contain an
automorphism of finite order ($\neq 1$). This contradicts Theorem 1
(cf. Remark 1.)
\end{proof}
If $\cF$ is a family of polynomials with coefficients in the field
$\overline{\Q}$ of all algebraic numbers, by $\cR_{\cF}$ we mean
the field obtained by adjoining to $\Bbb Q$ all the roots of the
polynomials in $\cF$.

Clearly Proposition 3 implies
\begin{corollary}
Let $\cF$ be a family of polynomials in
$\overline{\Q}[x]$ for which there exists a natural number $t$ such the
polynomials in $\cF$ have degree $\leq t$ and all coefficients of
the polynomials in $\cF$ have degree (with respect to $\Bbb Q$)
$\leq t$. Then $\Q(\cR_{\cF})$ is a proper subfield of the field $\overline{\Q}$ of
all algebraic numbers.
\end{corollary}
The next theorem will be proved in several steps.
\begin{theorem}\label{thm2}
For any proper subfield $K$ of $\QQ$ there exists an algebraic number $\alpha$ such that $\alpha^n$ does not lie in $K$ for any natural number $n$.
\end{theorem}
Theorem 2 is an immediate consequence of the following Theorem \ref{thm3},
whose proof depends on 
\begin{lemma}\label{lemma}
Let $\Q(\zeta_n)$ be the $n$-th cyclotomic field. There exists a number $\alpha \neq 0$ in  $\Q(\zeta_n)$ such that  for every automorphism $\sigma\in\emph{Gal}(\Q(\zeta_n)/\Q)$, $\sigma$ not the
identity, the quotient $\sigma(\alpha)/\alpha$ is not a root of unity.
\end{lemma}
\begin{proof}
(The following results from algebraic number theory can e.g. be
found in \cite{Marcus}.) There exist (infinitely many) prime numbers
$p\equiv 1 $ modulo $n$. In the ring of algebraic integers in $\Q(\zeta_n)$ any such prime $p$ splits into $\varphi(n) (= [\Q(\zeta_n):\Q])$ distinct
prime ideals of degree one. If $\mathfrak{p}$ is such a prime ideal, then
$\sigma(\mathfrak{p}) \neq \mathfrak{p}$ for every $\sigma\in Gal(\Q(\zeta_n)/\Q)$, $\sigma$ not the identity. Let $h$ be the class
number of $\Q(\zeta_n)$, then $\mathfrak{p}^h$ is a principal ideal. If
$\alpha$ is a generator of this ideal the conjugates of $\alpha$,
i.e., $\sigma(\alpha)$, $\sigma$ running through $Gal(\Q(\zeta_n)/\Q)$, generate distinct ideals. Hence for every $\sigma\in Gal(\Q(\zeta_n)/\Q)$, $\sigma$ not the identity, the quotient
$\sigma(\alpha)/\alpha$ is not a unit (i.e., an invertible element in the ring of
algebraic integers), in particular not a root of unity.
\end{proof}
\begin{theorem}\label{thm3}
Let $K$ be a field of characteristic 0 and
$L/K$ any extension where $K\subsetneq L$. Then there exists an
element $\alpha$ in $L$ such that no power $\alpha^n, n\geq 1$, lies in
$K$.
\end{theorem}

\begin{proof}
The assertion is trivial if $L$ is not algebraic over
$K$. So it suffices to consider the case where $L/K$ is a finite
algebraic extension.
We may obviously assume that $K$ is maximal subfield of
$L$, i.e., there is no field lying strictly between $K$ and
$L$. Furthermore we may  assume that there exists an element  $\beta
\in L\setminus  K$ for which  $\beta^n$ lies in $K$ for some integer $n > 1$.

We now distinguish  between two cases:

i) $L = K(\zeta)$ for some root of unity $\zeta$.

ii) $L \neq K(\zeta)$ for every root of unity $\zeta$.\\

\noindent
ad i) By Lemma 1 there exists an element $\alpha$ in $\Bbb Q(\zeta)$
such that for every $\sigma$ in Gal($\Bbb Q(\zeta)/\Bbb Q)$, $\sigma$ not the
identity, the quotient $\sigma(\alpha)/\alpha$ is not a root of unity.
This implies that $\sigma(\alpha^n) \neq \alpha^n$ for every $n\geq
1$. Therefore no power  $\alpha^n$, $n\geq 1$, lies in any proper subfield of $\Bbb
Q(\zeta)$. Clearly $ \Bbb Q(\alpha) =  \Bbb Q(\zeta)$ and hence $L = K(\alpha)$.

Since $K\cap \Bbb Q(\zeta)$ is a proper subfield of $\Bbb Q(\zeta)$
the element $\alpha$ constructed above has the property  that no power
$\alpha^n$, $n\geq 1$, lies in $K$.\\

\noindent
ad ii) By assumption there exists an element $\beta$ in $L\setminus K$
such  that $\beta^n$ lies in $K$ for some $n > 1$. We may assume that
$n$ is the smallest such number. The polynomial $f(x):= x^n - \beta^n$ is an
irreducible polynomial in $K[x]$. Indeed,  $f(x)= \prod_{j = 0}^{n-1}[x- \beta\zeta_n^j]$, where $\zeta_n$ is a primitive $n$-th root of unity. If $f(x)$ were reducible in $K[x]$, there would be an integer $t$, $1\leq t\leq
n-1$, such that for some $n$-th root of unity $\zeta$ the product $\beta^t\cdot \zeta$ would belong to $K$. Hence $ K(\beta^t) = K(\zeta)$. Since $K$ is a maximal
subfield of $L$ and $\beta^t$ is not in $K$, we would get $L =
K(\beta^t) = K(\zeta)$. But this contradicts our assumption ii). Hence
$f(x)$ is irreducible in $K[x]$ and thus the degree $[L:K]$ is $n$.

The above, in particular, shows that if some non-trivial power of an
element $\gamma\in L\setminus K$ lies in $K$, then $\gamma^n$ lies in
$K$ and $n$ is the smallest natural number with this property.

We apply this to  $\gamma : = 1+\beta$. No non-trivial power of
$\gamma$ lies in $K$. Otherwise, the above remark shows that
$\gamma^n$ would lie in $K$. But since
$\beta^n$ lies in $K$ the equation

$$\gamma^n = (1+\beta)^n = 1 + {n \choose 1} \beta +\dots + {n \choose n-1}\beta^{n-1}
+\beta^n$$

\noindent shows that  $\beta$ would be root of a non-zero polynomial in $K[x]$
of degree $n -1$, contradicting the fact that the 
degree of $\beta$ with respect to $K$ is $n$. Thus no power of $\gamma$
lies in $K$. 

The proof of Theorem 3 is now complete.
\end{proof}
From the above we deduce
\begin{corollary}\label{cor1}
If we view $\overline{\Q}$ as a subfield of
the complex number field $\C$, the algebraic numbers $\alpha$ such
that no power $\alpha^n$, $n\geq 1$, lies in a prescribed fixed
proper subfield $K$ of $\overline{\Q}$ form a dense subset of $\C$
(equipped with the usual topology).
\end{corollary}
\begin{proof}
 We distinguish between two cases: Either $K
\subseteq \Bbb R$ or $K\nsubseteq \Bbb R$.  

If  $K\subseteq \R$ then we notice that no power $\alpha^n$,
$n\geq 1$, of an integer $\alpha = a + b\cdot i $ in the Gaussian field $\Q(i)$
lies in $K$ if for instance $a$ is an odd rational
integer and $b$ is an even rational integer. (Just consider the prime
factorization of $\alpha$ and notice that $K\cap \Q(i) = \Q.)$
 Multiplying  all these numbers by rational numbers we get the desired subset of $\C$.

If $K\nsubseteq \R$ the numbers in $K$ lie dense in $\C$. So
if $\alpha$ is some algebraic number such that $\alpha^n\notin K$,
$n\geq 1$, then the numbers $\alpha\cdot k$, $k$ running through $K$,
yield the desired subset of $\C$.
\end{proof}
An immediate consequence of the previous results is
\begin{theorem}\label{thm3}
Let $n$ be a fixed natural number and $\cF$ the
family of all polynomials in $\Q[x]$ of degree $\leq n$.  Then
the algebraic numbers $\alpha$ such that no power $\alpha^n$, $n\geq 1$, lies in  $\Q(\cR_{\cF})$ form a dense subset of  $\C$.
\end{theorem}
Although not necessary for the following we point out another aspect
of Theorem 2. 
Let $\overline{\Q}^*$ be the multiplicative group of the non-zero 
elements of $\overline {\Q}$ and for a proper subfield $K$ of $\overline
{\Q}$ let $K^*$ be the multiplicative group of the  non-zero elements in
$K$.  Clearly the quotient group $\overline {\Q}^*/K^*$ is a divisible abelian
group. An additively written divisible abelian group is the direct sum
of copies of $\Q$ and copies of $\Bbb Z(p^{\infty})$, $p$ being a
prime and $\Z(p^{\infty})$ the divisible ($=$ injective) hull of
the cyclic group $C_p$ of order $p$. (See e.g. \cite{Fuc})

Theorem 1 states that for a proper subfield $K$ of $\QQ$ the quotient $\overline{\Q}^*/K^*$ is not a torsion group, in particular the corresponding (additively written)  divisible group contains at least
one copy of $\Q$. A sharper result is the following
\begin{theorem}
Let $K$ be a proper subfield of $\overline{\Bbb
Q}$. Then the torsion-free part of the divisible quotient group
$\overline{\Bbb Q}^*/K$ is isomorphic to the direct sum of
countable many copies of $\Bbb Q$. The torsion part of $\overline{\Bbb
Q}^*/K$ is non-trivial, i.e., contains at least one copy of $\Bbb
Z(p^{\infty})$ for some prime number $p$. 
\end{theorem}
\begin{proof}We first prove the assertion concerning the
torsion-free part of $\overline{\Bbb Q}^*/K.$

We start by considering the case where $K$ is real closed. For each
prime number $p_j \equiv 1$ (mod 4) let $\pi_j$ be an irreducible
factor of $p_j$ in the Gaussian ring $\Bbb Z[1,i]$. Then no product
of powers  ($\neq 1$) of  $\pi_j$'s  lies in $\Bbb Q$ and hence not in $K$
since $K\cap \Bbb Q(i) = \Bbb Q$.

Next assume $K$ is not real closed. Then $[\overline{\Bbb Q} : K]$ is
infinite; hence there exists an infinite tower of fields $K\subsetneq K_1\subsetneq K_2 \subsetneq \cdots \subsetneq K_i\subsetneq K_{i +1}\subsetneq \cdots $, where each field is a finite Galois extension of the previous.
\
By theorem 2 for each $i$ there exists an element $\alpha_i$ in
$K_{i+1}\setminus K_i$ such that no power ($\neq 1$) of $\alpha_i$
lies in $K_i$. Then no product of powers ($\neq 1)$ of the elements $\alpha_i, i\in \Bbb N,$ lies in $K$.  Therefore the corresponding
residue classes of the $\alpha$'s in $\overline{\Bbb Q}^*/K$ gives
rise to a direct sum of countably many of copies of $\Bbb Q$.

Concerning the torsion part of  $\overline{\Bbb Q}^*/K$ the assertion
is obvious if some root of unity does not lie in $K$. Hence we may
assume that all roots of unity belong to $K$.  Since, in particular,
$K$ is not real closed,  $\overline {\Bbb Q}/K$ is a proper infinite
Galois extension. This implies that there is an automorphism $\sigma\in
Gal(\overline {\Bbb Q}/K)$ of order $p^{\infty}$, for some prime number
$p$. Hence there is a Galois extension of $K$ whose Galois group is
the additive group of $\Bbb Z_p$ of $p$-adic integers (cf. \cite{Gey}). In particular,
there will be a cyclic extension $L/K$ of degree $p$. This extension
must be a Kummer extension $K(\root p\of \alpha)$ for some $\alpha\in
K$.
Thus  $\overline{\Bbb Q}^*/K$ contains an element of order $p$ and
therefore a copy of $\Bbb Z(p^{\infty})$.
\end{proof}

\section{Application to transcendental number theory}

Before giving an application of the previous results to transcendental
numbers we briefly recall some - by now classic - facts concerning
transcendental number theory.

At the 1900 International Congress of Mathematicians in Paris, Hilbert proposed his famous list of 23 problems  and the seventh of them asked about the arithmetic nature of the powers $x^y$, where both these numbers are algebraic. In 1934, Gelfond \cite{gel} and Schneider \cite{sch}, independently,
completely solved the problem.
\begin{theorem}[Gelfond-Schneider] \label{GelSchn}
Assume $\alpha$ and $\beta$ are algebraic numbers, with $\alpha \neq 0$ or $1$, and $\beta$ irrational. Then $\alpha^{\,\beta}$ is transcendental.
\end{theorem}
It follows that $2^{\sqrt{2}},\ \sqrt{2}^{\sqrt{3}}$ and $i^i$ are transcendental numbers. As well as $e^{\pi}$, because $e^{\pi}=(-1)^{-i}$. 

This theorem also classifies the arithmetic nature of powers of two algebraic numbers. However, there is no a similar result for powers $x^y$ when at least one of $x$ and $y$ is transcendental (see Table \ref{tab1}). In light of Theorem \ref{GelSchn} we might think that the power of two transcendental numbers is still transcendental, but this is not the case ($e^{\log 2}=2$).
\begin{table}[!htb]
\centering
\begin{tabular}{|c|c|c||p{3.7 cm}|}
\hline
\textbf{$x$} & \textbf{$y$} & \textbf{$x^y$} & Arithmetic nature\\
\hline \hline
2 & $\log 3/\log 2$ & 3 & Algebraic\\ \hline
2 & $i\log 3/\log 2$ & $3^i$ & Transcendental\\ \hline
$e^{i}$ & $\pi$ & -1 & Algebraic\\ \hline
$e$ & $\pi$ & $e^{\pi}$ & Transcendental\\ \hline
$2^{\sqrt{2}}$ & $\sqrt{2}$ & 4 & Algebraic\\ \hline
$2^{\sqrt{2}}$ & $i\sqrt{2}$ & $4^{i}$ & Transcendental \\ \hline
\end{tabular}
\caption{Possibilities}
\label{tab1}
\end{table}\\
The case $x=y$ seems to be more interesting: can the number $T^T$ be algebraic for some transcendental $T$? Sondow and the second author \cite{SM} showed that the answer for this question is yes, actually they proved that
\newpage
\begin{proposition}[\textit{Cf. Proposition 1 in \cite{SM}}]
Given $A\in [e^{-1/e},\infty)$, let $T\in\mathbb{R}^+$ satisfy $T^T=A$. If either

\emph{(i).} $A^n\in \QQ \setminus\mathbb{Q}$ for all $n\in \N$, or

\emph{(ii).} $A\in \Q\setminus\{n^n:n\in \N\}$,\\
then $T$ is transcendental. In particular, $T\notin \QQ$ if $T^T\in\mathbb{Q}\cap (e^{-1/e},1)$.
\end{proposition}
However, for instance, they were not able to prove the existence of algebraic numbers which can be written of the form $T^{T^{2009}+T+1}$, with $T$ transcendental. Now using some results from the previous section we are able to solve completely this kind of problem.
\begin{theorem}\label{PQ}
For arbitrary non-constant rational functions $P(x),
Q(x)\in \overline{\Bbb Q}(x)$ the set of algebraic numbers of the form
$P(T)^{Q(T)}$ with $T$ transcendental, is  dense in some open subset
of the complex plane.
\end{theorem}
\begin{proof}
The set of complex numbers for which $P(x)$ or $Q(x)$ has
a pole or zero or $P(x)$ takes the value 1 is finite. The complement
of this set inside $\Bbb C$ is an open subset of the complex plane.
Let $\Omega$ be an open simply connected subset of the above open set.
Choosing, for instance, the principal branch of the multi-valued
logarithm  function, the function  $ f(x): = P(x)^{Q(x)}$ is well defined and
analytical in $\Omega$. This function is a non-constant function. Indeed, since
$P(x)$ and $Q(x)$ are non-constant and the algebraic numbers form a
dense subset of $\Bbb C$ there exists an algebraic number $\beta$ in
$\Omega$ for which $Q'(\beta)P(\beta) \neq 0$. If $f(x)$ were
constant then $\log P(x) = - Q(x)P'(x)/Q'(x)P(x)$. Setting $x = \beta$
we get that $\log P(\beta)$ would be an algebraic number. But this would
contradict the famous theorem by Lindemann (cf [6]) that log $\beta \notin
\overline{\Bbb Q}$ for all $\beta\in \overline{\Bbb Q}\setminus
\{0,1\}$.
Since a non-constant analytic function maps an open connected set of $\Bbb C$ onto an open
connected set of $\Bbb C$ we see that $f(\Omega)$ is an open connected
subset of $\Bbb C$.

Write $Q(x)$ as $Q_1(x)/Q_2(x)$, where $Q_1(x)$ and $Q_2(x)$ are
polynomials in $\overline{\Bbb Q}[x]$. Let $\cF$ be the family of
polynomials of the form $Q_1(x) - dQ_2(x)$, $d$ running through $\Bbb
Q$.  Clearly there exists an integer $t$ such that every polynomial in
$\cF$ has degree $\leq t$ and all coefficients of these
polynomials have a degree (with respect to $\Bbb Q$) $\leq t$.

Hence Corollary 1 implies that $\Q(\cR_{\cF}) \neq \overline{\Bbb
Q}$. By Corollary 2 the algebraic numbers $\alpha$ such that no power
$\alpha^n$, $n\geq 1$ lies in $\Q(\cR_{\cF})$ form a dense subset
of $\Bbb C$. Since $f(\Omega)$ is open, it contains dense subset of
numbers $\alpha$ with the above property. Every $\alpha$
in this dense set has the form $\alpha = f(T) = P(T)^{Q(T)}$. This
number $T$ must be transcendental. In fact, assume $T$ were algebraic.
Since $P(T)\notin \{0,1\}$ by the Gelfond-Schneider theorem we
conclude that $Q(T)$ would be a rational number say $d = \frac{r}{s}$,
$r$ and $s$ integers and $s >0$.  
Hence $T$ would  belong to $\cR_{Q_1(x)-\frac{r}{s}Q_2(x)} \subseteq \cR_{\cF}$. But then $\alpha^s = P(T)^r$ and thus $\alpha^s \in \Q(\cR_{\cF})$ contradicting the choice of $\alpha$.
\end{proof}
\begin{remark}
Note that for any transcendental number $T$ and any non-constant rational function $P\in \QQ(x)$, the number $P(T)$ is transcendental. Otherwise, if we write $P = P_1/P_2$, where $P_1$ and $P_2$ are polynomials in $\overline{\Q}[x]$ the number $T$ would be a root of the polynomial $P_1(x) -AP_2(x) \in \overline{\Q}[x]$. But any such root is algebraic since $\overline{\Q}$ is algebraically closed.
\end{remark}
\begin{exe}
If $P(x)=x$ and $Q(x)=x^2$, then easy calculations show that the set of the algebraic numbers of the form $T^{T^2}$, with $T$ transcendental, is dense in the interval $[1/\sqrt{e},+\infty)$.
\end{exe}
In light of what we just proved, the following problem arises: is there algebraic number in the form $T^{T^T}$, for some transcendental $T$? And in the form $T^{T^{T^T}}$? And so on? 

This question is still open. The main open problem in transcendental number theory is the renowned
Schanuel's conjecture (see \cite{Lang} for its statement). In his Ph.D thesis \cite[p. 32]{Mar}, the second author proved that
\begin{proposition}
Suppose that Schanuel's conjecture is true. For any $m\geq 3$ and any algebraic number $1 < A\notin \N$, there exists a transcendental number $T$ such that $A=\underbrace{T^{\,T^{\,\cdot^{\cdot^{\cdot^{T}}}}}}_{m}$.
\end{proposition}
We finish by a related question which may be considered as an inverse
problem of Theorem \ref{PQ}.
\begin{question}
Give example, if any, of a well-known transcendental number $T$ (like $e, \pi, \log 2$, etc) such that there exist non-constant rational functions $P(T),Q(T)\in \QQ(x)$, such that $P(T)^{Q(T)}$ is algebraic.
\end{question}
In this connection, the second author \cite{Mar2} proved that if Schanuel's conjecture is true then $P(e)^{Q(e)}$ is transcendental, for any non-constant polynomials $P(T),Q(T)\in \QQ[x]$.

\section*{Acknowledgement}
\noindent
Part of this work was started when the second author was still a doctorate student at Universidade de Bras\' ilia. He thanks this institution for its hospitality.

\end{document}